\documentstyle[12pt]{article}
\begin{document}
\title{ Some extensions of Hardy's integral inequalities to  Hardy type
spaces
}
\author{{    Shunchao Long }\\
 }
\date{}
\maketitle
\begin{center}
\begin{minipage}{120mm}
\vskip 0.1in

{
\begin{center}
{\bf Abstract}
\end{center}
~~  In this paper some extensions of Hardy's integral inequalities
to $0<p\leq 1$ are established.}
\end{minipage}
\end{center}
\vskip 0.1in
\baselineskip 16pt
\par \section*
 {\bf  1.  Introduction}

 \par Let
$$Hf(x)=\frac {1}{x}\int^x_0f(t)dt, ( x>0), $$
and (dual form)
$$H^*f(x)=\int^{\infty}_xf(t)dt, (x>0).$$
\par The   Hardy integral inequality can be stated as (see [3]):
 \begin{equation}
  \int^{\infty}_0(Hf(x))^pdx< p'^p\int^{\infty}_0(f(x))^pdx,  p>1, f(x)\geq 0 ,
\end{equation}
\begin{equation}
 \int^{\infty}_0(H^*f(x))^pdx< p^p\int^{\infty}_0(xf(x))^pdx,  p>1, f(x)\geq 0,
\end{equation}
 unless $f\equiv 0.$
\par The inequality (1)  was firstly proved by Hardy [4], but the constant was not
determined; and Landau found out precisely the constant  is $p'^p$
in [7]; later,
 Hardy [5] generalized it to the inequality (2) himself. However, if $0<p< 1$, the reverse direction
 inequalities  hold (see [3]):
 \begin{equation}
  \int^{\infty}_0(Hf(x))^pdx> p'^p\int^{\infty}_0(f(x))^pdx,  0<p< 1, f(x)\geq 0,
\end{equation}
\begin{equation}
 \int^{\infty}_0(H^*f(x))^pdx> p^p\int^{\infty}_0(xf(x))^pdx, 0<p< 1, f(x)\geq 0 .
\end{equation}
 unless $f\equiv 0.$
 The constants in (1),(2), (3) and (4) are the best  possible.

\par The positive direction inequalities (1) and (2) play an
important role in many areas such as  harmonic analysis [12], PDE
[8], etc.
 In view of this,  much efforts and
  time have been devoted to their improvement and generalizations over the years.
\par When $p>1$, many generalizations  include the works in  numerous papers, for example,
   [3,11,6] and some of the references cited therein.
\par When $p=1$, in view of the theory of Hardy spaces on ${\bf R}^n$ established by
Coifman and Weiss in [1] and some others, J. Garcia-Cuerva and J. L.
Rubio de Francia extended (1) to Hardy spaces, and established a
positive direction inequality of Hardy type for $p=1$: let $f\in
H^1({\bf R})$ be supported in $[0,\infty)$, then
\begin{equation}
  \int^{\infty}_0|Hf(x)| dx<({\rm log~} 2
  )~\|f\|_{H^{1,\infty}_{at}({\bf R})},
\end{equation}
where $H^{1,\infty}_{at}$ is the atom Hardy spaces. See [2].

 In this paper we extend the positive direction inequalities (1), (2) and (5) to $0<p\leq 1$, as well as establish
 some  estimates of $H$ and $H^*$ from Hardy type spaces to Hardy type spaces.

Let us introduce some definitions of the  Hardy type spaces.
\par {\bf Definition 1} ~~Let $ 0<p\leq  1 \leq q \leq \infty, p<q $
and $ s \in {\bf N} $ and $  w\ge 0 $ is a weight function on $ {\bf
R}^+$.
\par (a) A function $a(x)$ on $   {\bf R}^+ $ is said to be a   $ (p ,q , s)_w$-atom , if
\par (i)~~~~supp $a\subset  (x_0, x_1) \subset  {\bf R}^+, x_0>0, $
\par (ii)~~~~$\|a\|_{L^{q}_w({\bf R}^+ )}\leq \left( \int ^{x_1}_{x_0}w(x)dx\right)^{1/q-1/p},$
\par  (iii)~~~~  $\int_{{\bf R^+}}a(x)x^\beta dx=0, \beta =0,1, ..., s,$
 \par  (b) and $a(x)$ is said to be a  $L-  (p ,q ,s )_w$-atom , if it is  a $ (p ,q , s)_w$-atom and satisfies
\par  (iv)~~~~$\int_{{\bf R^+}}a(x){\rm ln} x dx=0.$
 \par {\bf Definition 2}~~ Let $ p, q , s  $ and $  w  $ as in Definition 1.
 Some Hardy spaces on ${\bf R^+}$ are defined by
\par
$H ^ {p,q,s}_w({\bf R^+} )=\{ f: f=\sum _{k=1}^{\infty} \lambda
_ka_k $, where each $a_k$ is a $ (p ,q , s)_w$-atom, $\sum
_{k=1}^{\infty} |\lambda _k|^p <+ \infty $, and the series converges
in the sense of distributions\},
\par and
\par
$LH ^ {p,q,s}_w({\bf R^+} )=\{ f: f=\sum _{k=1}^{\infty} \lambda
_ka_k $, where each $a_k$ is a $L- (p ,q , s)_w$-atom, $\sum
_{k=1}^{\infty} |\lambda _k|^p <+ \infty $, and the series converges
in the sense of distributions\}.
 \par And define the quasinorms of a function of $H ^ {p,q,s}_w({\bf R^+} )$
 or
$LH ^ {p,q,s}_w({\bf R^+} )$ by
\par
 $\|f\|_{H ^{ p,q,s}_{w}({\bf R^+} )}= \inf \left(\sum\limits_{k=1}^{\infty}|\lambda _k|^p\right)^{1/p},$  ~~or~~
$\|f\|_{LH ^{ p,q,s}_{w}({\bf R^+} )}= \inf
\left(\sum\limits_{k=1}^{\infty}|\lambda _k|^p\right)^{1/p} $
\\
respectively, where the infimum is taken over all the decompositions
of $f$ as above.

\par Simply, we  denote   $H _1^ {p,q,s} ({\bf R^+} )$ and $LH
_1^ {p,q,s} ({\bf R^+} )$ by $H ^ {p,q,s} ({\bf R^+} )$ and $LH ^
{p,q,s} ({\bf R^+} )$ respectively.
\par $H ^ {p,q,0} ({\bf R^+} )$ is an example of the Hardy spaces on spaces of homogeneous type  studied
 by Coifman and Weiss [1], and A.Marcias and C. Segovia
 [9,10].

\par Throughout this paper, we denote $p'=p/(p-1)$ for $1<p<\infty, \infty '=1,$ and $ 1'=\infty,$ and suppose
all functions are not identical to zero.

\par We extend (1) and (5) to the case of $0<p\leq 1$:
 \par  {\bf Theorem 1}~~Let $ 0<p\leq  1 \leq q \leq \infty$ and $p<q.  $ If $f$ is in
$H ^ {p,q,0} ({\bf R^+} )$,  then,
\begin{eqnarray}
\|Hf \|_{L^p({\bf R}^+)}< \frac{1}{(1-p/q)^{1/p}}\|f \|_{H ^ {p,q,0}
({\bf R^+} )}.
\end{eqnarray}

\par We extend (2)  as well:
\par  {\bf Theorem 2}~~
 Let $ 0<p\leq  1  .  $ If $f$ is in
$H ^ {p,q,0}_{x^p} ({\bf R^+} )$,  then,
 \begin{eqnarray}
 \|Hf \|_{L^p({\bf R}^+)}<
\left\{\begin{array}{cc}
 \|f \|_{H ^ {p,q,0}_{x^p} ({\bf R^+} )} , &{\rm if } ~~q=\infty,\\
\frac{1}{(1-p )^{1/p}(p+1) } \|f \|_{H ^ {p,q,0}_{x^p} ({\bf R^+}
)},
&{\rm if }~~ q=1, p\not= 1,\\
\frac{1}{\left(1 -pq'/q \right)^{1/q'}}
 \frac{ ( {1+p} )^{ 1/q-1/p}
}{((1/q'-p/q)p+1)^{1/p}} \|f \|_{H ^ {p,q,0}_{x^p} ({\bf R^+} )},
&{\rm if }~~ 1<q<\infty,p<q-1.
\end{array}
\right.
\end{eqnarray}

\par We also extend these results to the estimates of $H$ and $H^*$
 from Hardy type spaces to Hardy type spaces.
 \par  {\bf Theorem 3}~~Let $ 0<p\leq  1 < q \leq \infty$ and $ s \in {\bf N} $. If $ f$ is in
$H ^ {p,q,s} ({\bf R^+} )$, then,
\begin{eqnarray}
\|Hf\|_{H ^ {p,q,s}({\bf R^+})}\leq q'\|f\|_{LH
^ {p,q,s}({\bf R^+})}.
\end{eqnarray}
 \par  {\bf Theorem 4}~~Let $ 0<p\leq  1 \leq q \leq \infty,  s \in {\bf N} $ and $s-1>0$. If $  f$ is in
$H ^ {p,q,s} ({\bf R^+} )$, then,
\begin{eqnarray}
 \|H^*f\|_{H ^ {p,q,s-1}({\bf R^+})}\leq
\left\{\begin{array}{cc}
 (1+p)^{1/p}\|f\|_{H
^ {p,q,s}({\bf R^+})} , &{\rm if } ~~q=\infty,\\
\frac{1}{(1-p)(1+p)^{1-1/p}}\|f\|_{H
^ {p,q,s}({\bf R^+})}, &{\rm if }~~ q=1, p\not= 1,\\
q\|f\|_{H
^ {p,q,s}({\bf R^+})}, &{\rm if }~~ 1<q<\infty,0<p\leq 1.
\end{array}
\right.
\end{eqnarray}

\par \section*{
\bf  2.  Proof of Theorems}
 \par~~~~ Firstly, let us introduce two inequalities.   By Minkowski inequalities [3],   it is easy to see that
\begin{eqnarray}
{(  {x_1 }-{x_0}  )^{p+1}}< {x_1^{p+1}-x_0^{p+1}};
\end{eqnarray}
 when $x_1>x_0>0$ and $ p > 0,$ and
\begin{eqnarray}
{(  {x_1 }-{x_0}  )^{1-p}}> {x_1^{1-p}-x_0^{1-p}}.
\end{eqnarray}
when $x_1>x_0>0$ and $1> p > 0.$
\par  {\bf Proof of Theorem 1}~~It suffices to  prove
  the following propositions.
\par {\bf Proposition 1} ~~Let $ 0<p\leq  1 \leq q \leq \infty$ and
$p<q.$  We have
\begin{equation}
 \int^{\infty}_0|Ha(x)|^pdx< \frac{1}{1-p/q}
\end{equation}
for all $ (p,q,0)$-atom $a$ on ${\bf R}^+$.

\par Let $FH ^ {p,q,0} ({\bf R^+} )$ be the set of all finite linear
combination of $ (p,q,0)$-atoms.

\par From Proposition 1 it is easy to get that
\par {\bf Proposition 2} ~~Let $ 0<p\leq  1 \leq q \leq \infty$ and
$p<q.$ Then (6) holds for $f$ in $FH ^ {p,q,0} ({\bf R^+} )$.

\par {\bf Proposition 3} ~~Let $ 0<p\leq  1 \leq q \leq \infty$ and
$p<q.$  Then
$$Hf(x)=\sum _{k=1}^{\infty} \lambda _kHa_k(x) ~~~~{\rm a.e.}$$
for all $f=\sum _{k=1}^{\infty} \lambda _ka_k \in H ^ {p,q,0} ({\bf
R^+} )$ , where each $a_k$ is a $  (p ,q , 0) $-atom and $ \sum
_{k=1}^{\infty} |\lambda _k|^p<\infty$.

In fact, once Proposition 1 and Proposition 3 have been proved,
then, for  $f  \in
 H ^ {p,q,0} ({\bf R^+} )$, i.e. $f=\sum _{k=1}^{\infty}
\lambda _ka_k $, we have
\begin{eqnarray*}
\int^{\infty}_0|Hf(x)|^pdx &= &
\int^{\infty}_0|\sum\limits_{k=1}^{\infty} \lambda _kHa_k (x)|^pdx
\\
&\leq & \int^{\infty}_0 \sum\limits_{k=1}^{\infty} |\lambda _k |^p |
Ha_k (x)|^pdx
 \\
&\leq & \sum\limits_{k=1}^{\infty} |\lambda _k |^p  \int^{\infty}_0|
Ha_k (x)|^pdx
\\
& & {\rm (by~~ Minkowski~~ inequality)}
\\
&< & \frac{1}{1-p/q}\sum\limits_{k=1}^{\infty} |\lambda _k |^p ,
\end{eqnarray*}
   From this, (6) follows easily. Thus, we finish the proof of Theorem
1.

\par  {\bf Proof of Proposition 1}~~
Let $a$ be a   $ (p,q,0)$-atom on ${\bf R}^+$, i.e.  supp $a\subset
(x_0, x_1)\subset  {\bf R}^+, (x_0>0) $, $\|a\|_{L^{q} ({\bf R}^+
)}\leq \left(  {x_1}-{x_0} \right)^{1/q-1/p},$ and $\int _0^\infty
a(x)dx =0.$ Let us prove (12).
\par From the vanishing property of $a$, it is easy to see that  supp $Ha\subset (x_0, x_1) $.
Thus, noticing that $x_0>0,$     we have
\begin{eqnarray*}
\int^{\infty}_0|Ha(x)|^pdx &=& \int^{x_1}_{x_0}|Ha(x)|^pdx
\\
&=&
\int^{x_1}_{x_0}\left|\frac{1}{x}\int^{x }_{x_0} a(t)dt\right|^pdx
\\
&\leq &
\int^{x_1}_{x_0}\left(\frac{1}{x}\left(\int^{x }_{x_0} |a(t)|^qdt\right)^{1/q} (  {x }-{x_0}  )^{1/q'}\right)^pdx
\\
& &~~~~~~~~~~~~~~~~~~~{\rm (by~~ Holder's~~ inequality)}
\\
&\leq &
\|a\|^p_{L^{q} ({\bf R}^+ )}\int^{x_1}_{x_0}\left(\frac{1}{x}    (  {x }-{x_0}  )^{1/q'}\right)^pdx
\\
&\leq &
({x_1}-{x_0}  )^{p/q-1 }\int^{x_1}_{x_0} \frac{(  {x }-{x_0}  )^{p}}{x^p}({x }-{x_0}  )^{-p/q} dx
\\
&< &
({x_1}-{x_0}  )^{p/q-1 }\int^{x_1}_{x_0}  ({x }-{x_0}  )^{-p/q} dx
\\
& = &  1/(1-p/q)
\end{eqnarray*}
when $1< q <\infty;$
\begin{eqnarray*}
\int^{\infty}_0|Ha(x)|^pdx  &=&
\int^{x_1}_{x_0}\left|\frac{1}{x}\int^{x }_{x_0} a(t)dt\right|^pdx
\\
&\leq &
\|a\|^p_{L^{q} ({\bf R}^+ )}\int^{x_1}_{x_0}\left(\frac{1}{x} \right)^pdx
\\
&\leq &
 ({x_1}-{x_0}  )^{p -1 }\int^{x_1}_{x_0} \frac{(  {x }-{x_0}  )^{p}}{x^p}({x }-{x_0}  )^{-p } dx
\\
&< &
({x_1}-{x_0}  )^{p -1 }\int^{x_1}_{x_0}  ({x }-{x_0}  )^{-p } dx
\\
&= & 1/(1-p )
\end{eqnarray*}
when $   q =1$ and $p\neq 1$; and

\begin{eqnarray*}
\int^{\infty}_0|Ha(x)|^pdx
&=&
\int^{x_1}_{x_0}\left|\frac{1}{x}\int^{x }_{x_0} a(t)dt\right|^pdx
 \\
&\leq &
\|a\|^p_{L^{\infty} ({\bf R}^+ )}\int^{x_1}_{x_0}\left(\frac{1}{x}    (  {x }-{x_0}  )^{ }\right)^pdx
\\
&\leq &
({x_1}-{x_0}  )^{ -1 }\int^{x_1}_{x_0} \frac{(  {x }-{x_0}  )^{p}}{x^p} dx
\\
&< &
  1 .
\end{eqnarray*}
when $  q =\infty.$
  Thus,(12) has been proved for $1\leq q \leq\infty,$ and $ p<q$.
  Thus, we finish the proof of Proposition 1.

  \par  {\bf Proof of Proposition 3}~~ Let $H ^ {p,q,0} ({\bf
R^+} )$, then $f=\sum _i \lambda _i a_i
 $ where $a_i$ are  $(p,q,0) $-atoms and $\sum_i
 |\lambda_i|^{p} <+\infty.$ We know that $H a_i$ is well defined
 for every $i$ since $a_i\in L^q({\bf
R^+} )$ with $1\leq q< \infty$ and (12) holds
  by Proposition 1.  Then
$$
\|\sum _i \lambda _i H a_i\|^{p}_{L^{p }({\bf R^+} )}
 \leq   \sum _i |\lambda _i |^{p} \|H a_i\|^{p}_{L^{p }({\bf
R^+} )}
  \leq C \sum _i |\lambda _i |^{p} \leq C \|f\|_{ {H}^{p,q,0}({\bf
R^+} )}^{p}<\infty,
  $$
it follows $|\sum _i \lambda _i H a_i(x)|<\infty ~ {\rm
  a.e. }$.
 Let
\begin{eqnarray}
f=\sum _i \lambda _i^{(1)} a_i^{(1)}=\sum _i \lambda _i^{(2)}
a_i^{(2)}
\end{eqnarray}
 with
\begin{eqnarray}
\sum_i
 |\lambda_i^{(1)}|^{\bar{p}} <+\infty
 {~~\rm and ~~} \sum_i
 |\lambda_i^{(2)}|^{\bar{p}} <+\infty
 \end{eqnarray}
 and  $a_i^{(1)}$ and  $a_i^{(2)}$ are  $(p,s,\alpha) $-atoms.
Once it is proved that
\begin{eqnarray*}
 \sum _i \lambda _i^{(1)} Ha_i^{(1)}=\sum _i \lambda _i^{(2)}
Ha_i^{(2)}  ~~{\rm a.e.},
 \end{eqnarray*}
then,
  $$H f(x) =\sum _i \lambda _i H a_i(x)
~  {\rm
  a.e. }$$
 is well defined  for  all $f=\sum _{k=1}^{\infty} \lambda _ka_k \in H ^ {p,q,0} ({\bf
R^+} )$. Thus, Proposition 2 holds.

It is remained to prove (13). For any $\delta>0$, by (14), there
exists $i_0$ such that
\begin{eqnarray}
\sum_{i=i_0}^\infty
 |\lambda_i^{(1)}|^{\bar{p}} <\delta^{\bar{p}}
 {~~\rm and ~~} \sum_{i=i_0}^\infty
 |\lambda_i^{(2)}|^{\bar{p}} <\delta^{\bar{p}}.
 \end{eqnarray}
From (13), we see that
$$
\sum _{i=1}^{i_0-1} (\lambda _i^{(1)} a_i^{(1)}- \lambda _i^{(2)}
a_i^{(2)} )= \sum _{i=i_0}^{\infty}\lambda _i^{(2)} a_i^{(2)}-\sum
_{i=i_0}^{\infty} \lambda _i^{(1)} a_i^{(1)},
 $$
then,
\begin{eqnarray}
\|\sum _{i=1}^{i_0-1} (\lambda _i^{(1)} a_i^{(1)}- \lambda _i^{(2)}
a_i^{(2)} )\|^{p}_{ H^{p,q,0}({\bf R^+} )} \leq \sum_{i=i_0}^\infty
 |\lambda_i^{(1)}|^{p} + \sum_{i=i_0}^\infty
 |\lambda_i^{(2)}|^{p} <2 \delta^{p}.
 \end{eqnarray}
By the linearity of $H$, we have
\begin{eqnarray*}
 \sum _{i=1}^{\infty} \lambda _i^{(1)} H a_i^{(1)}- \sum
_{i=1}^{\infty} \lambda _i^{(2)} H a_i^{(2)}
  =H (\sum
_{i=1}^{i_0-1} (\lambda _i^{(1)} a_i^{(1)}- \lambda _i^{(2)}
a_i^{(2)} ))+ \sum _{i=i_0}^{\infty}\lambda _i^{(1)} H
a_i^{(1)}-\sum _{i=i_0}^{\infty} \lambda _i^{(2)} H a_i^{(2)}.
\end{eqnarray*}
\begin{eqnarray}
\end{eqnarray}
By Proposition 2 and (16), we see that
\begin{eqnarray}
\|H (\sum _{i=1}^{i_0-1} (\lambda _i^{(1)} a_i^{(1)}- \lambda
_i^{(2)} a_i^{(2)} ))\|^{p}_{{L}^{p}({\bf R^+} )} \leq \|\sum
_{i=1}^{i_0-1} (\lambda _i^{(1)} a_i^{(1)}- \lambda _i^{(2)}
a_i^{(2)} )\|^{p}_{ {H}^{p,s}({\bf R^+} )}<2 \delta^{p}.
\end{eqnarray}
From (17), (18), (12) and  (15), we have
$$
\|\sum _{i=1}^{\infty} \lambda _i^{(1)} T_\varepsilon a_i^{(1)}-
\sum _{i=1}^{\infty} \lambda _i^{(2)} T_\varepsilon a_i^{(2)}
\|^{p}_{{L}^{p}({\bf R^+} )} <4 \delta^{p}.
 $$
Let $\delta\rightarrow 0$, we get that $ \|\sum _{i=1}^{\infty}
\lambda _i^{(1)} H a_i^{(1)}- \sum _{i=1}^{\infty} \lambda _i^{(2)}
H a_i^{(2)} \|^{p}_{{L}^{p}({\bf R^+} )} =0,
 $ it follows that $
\sum _{i=1}^{\infty} \lambda _i^{(1)} H a_i^{(1)}=\sum
_{i=1}^{\infty} \lambda _i^{(2)} H a_i^{(2)}  {\rm
  a.e. }.
 $ Thus, we finish the proof of Proposition 3.
\par  The proof of Theorem 1 is finished.
\par  {\bf Proof of Theorem 2}~~ As the proof of Theorem 1, it suffices to  prove
the following  propositions.

\par {\bf Proposition 4} ~~ Let $ 0<p\leq  1. $ We have
\begin{eqnarray}
 \int^{\infty}_0|H^*a(x)|^pdx<
\left\{\begin{array}{cc}
1, &{\rm if }~~ q=\infty,\\
\frac{1}{(1-p )(p+1)^{ p}}, &{\rm if }~~ q=1, p\not= 1,\\
\frac{1}{\left(1 -pq'/q \right)^{p/q'}}
 \frac{ ( {1+p} )^{ p/q-1}
}{(1/q'-p/q)p+1} , &{\rm if }~~ 1<q<\infty,p<q-1,
\end{array}
\right.
\end{eqnarray}
 for all $ (p,q,0)_{x^p}$-atom on ${\bf R}^+$.

\par Let $FH ^ {p,q,0}_{x^p}
({\bf R^+} )$ be the set of all finite linear combination of $
(p,q,0)_{x^p}$-atoms.

\par From Proposition 4 it is easy to get that

\par {\bf Proposition 5} Let $p$ and $q$ as in Theorem 2. Then (7)
holds for all $f\in FH ^ {p,q,0}_{x^p} ({\bf R^+} )$.

\par {\bf Proposition 6} ~~Let $ 0<p\leq  1, 1\leq q \leq \infty, p<q-1, p\neq q=1   $ and
$p<q.$  Then
$$Hf(x)=\sum _{k=1}^{\infty} \lambda _kHa_k(x) ~~~~{\rm a.e.}$$
for all $f=\sum _{k=1}^{\infty} \lambda _ka_k \in  H ^ {p,q,0}_{x^p}
({\bf R^+} )$, where each $a_k$ is a $ (p,q,0)_{x^p}$-atom on ${\bf
R}^+$ and $ \sum _{k=1}^{\infty} |\lambda _k|^p<\infty$.

\par  {\bf Proof of Proposition 4}~~
Let $a$ be a   $ (p,q,0)_{x^p}$-atom on ${\bf R}^+$ with supp
$a\subset (x_0, x_1)\subset  {\bf R}^+, (x_0>0) $. Let us prove
(19).

 \par As those discussions  in the proof of Theorem 1 we have   supp $H^*a\subset (x_0, x_1) $. Thus,
when $  q =\infty,$  noticing  that  $\|a\| _{L^\infty ({\bf
R}^+)}=\|a\| _{L^\infty_{t^p}({\bf R}^+)}, $   by (10), we have
\begin{eqnarray*}
\int^{\infty}_0|H^*a(x)|^pdx
&=& \int^{x_1}_{x_0}|H^*a(x)|^pdx
\\
&=&
\int^{x_1}_{x_0}\left| \int^{x_1 }_{x} a(t)dt\right|^pdx
 \\
 &\leq &
\|a\|^p_{L^\infty_{t^p}({\bf R}^+)}
  \int^{x_1}_{x_0}    ({x_1}-{x }  )^{ p }
  dx
\\
 &\leq &
 \left(\int^{x_1 }_{x_0}  t^{p }dt
\right)^{-1}
 \int^{x_1}_{x_0}    ({x_1}-{x }  )^{ p }
  dx
\\
&\leq &
   \frac{(  {x_1 }-{x_0}  )^{p+1}}{x_1^{p+1}-x_0^{p+1}}
\\
&< &
1;
\end{eqnarray*}
when $  q =1, 0<p<1$,  by (11), we have
\begin{eqnarray*}
\int^{\infty}_0|H^*a(x)|^pdx
  &\leq &
\int^{x_1}_{x_0}\left| \int^{x_1 }_{x} |a(t)|t^p \frac{1}{t^p}dt\right|^pdx
 \\
 &\leq &
\|a\|^p_{L^1_{t^p}({\bf R}^+)}
  \int^{x_1}_{x_0}    \frac{1}{x^{p^2}}
  dx
\\
 &\leq &
 \left(\int^{x_1 }_{x_0}  t^{p }dt
\right)^{p-1}
 \int^{x_1}_{x_0}    \frac{1}{x^{p^2}}
  dx
\\
&= &
\frac{(1+p)^{1-p}}{(1-p^2)} \frac{  x_1 ^{1-p^2}-x_0^{1-p^2}  }
{\left(x_1^{1+p}-x_0^{1+p}\right)^{1-p}}
\\
&< &
\frac{1}{(1-p )(1+p)^{ p}} \frac{ x_1 ^{1-p^2}-x_0^{1-p^2} }{ x_1 ^{1-p^2}-x_0^{1-p^2} }
\\
&= &
\frac{1}{(1-p )(1+p)^{ p}};
\end{eqnarray*}
  and when $1< q <\infty,$      $0<p<q-1,  $  we see that  $(p+1)/q<1,$ and then  $pq'/q<1,$
  by (10) and (11), we have
 \begin{eqnarray*}
\int^{\infty}_0|H^*a(x)|^pdx
&=&
\int^{x_1}_{x_0}\left| \int^{x_1 }_{x} a(t)dt\right|^pdx
\\
&\leq &
\int^{x_1}_{x_0}
\left(
 \left(\int^{x_1 }_{x} |a(t)|^qt^pdt\right)^{1/q}
\left(\int^{x_1 }_{x} \frac{1}{t^{pq'/q}}dt\right)^{1/q'}
\right)^pdx
\\
& &~~~~~~~~~~~~~~~~~~~{\rm (by~~ Holder's~~ inequality)}
\\
&\leq &
\|a\|^p_{L^q_{t^p}({\bf R}^+)}
\left(\frac{1}{1-pq'/q}\right)^{p/q'}
\int^{x_1}_{x_0}
\left(x_1^{1-pq'/q}-x ^{1-pq'/q}\right)^{p/q'}
  dx
 \\
&\leq &
\left(\frac{1}{1-pq'/q}\right)^{p/q'}
\left(\int^{x_1 }_{x_0}  t^{p }dt
\right)^{1-p/q}
 \int^{x_1}_{x_0}
\left((x_1^{ }-x) ^{1-pq'/q}\right)^{p/q'}
  dx
\\
&= &
  \left(\frac{1}{1 -pq'/q}\right)^{p/q'}
 \frac{ ( {1+p} )^{ p/q-1}
}{(1/q'-p/q)p+1}
 \frac{(  {x_1 }-{x_0}  )^{(1/q'-p/q)p+1}}{\left(x_1^{p+1}-x_0^{p+1}\right)^{1-p/q}}
 \\
&< &
  \left(\frac{1}{1 -pq'/q}\right)^{p/q'}
 \frac{ ( {1+p} )^{ p/q-1}
}{(1/q'-p/q)p+1}
 \frac{(  {x_1 }-{x_0}  )^{(1/q'-p/q)p+1}}{\left(x_1^{ }-x_0^{ }\right)^{(p+1)(1-p/q)}}
\\
&= &
   \frac{1}{\left(1 -pq'/q \right)^{p/q'}}
 \frac{ ( {1+p} )^{ p/q-1}
}{(1/q'-p/q)p+1}
  .
\end{eqnarray*}
   Thus,(19) has been proved for $1\leq q \leq\infty,$ and $ p<q-1$.

\par     Proposition 5 follows easily from Proposition 4.
Proof of Proposition 6 is same as that of Proposition 3.

\par  Thus, we finish the proof of Theorem 2.

\par {\bf Proof of Theorem 3}~~ Let $a$ be a   $ L-(p,q,s)$-atom on ${\bf R}^+$  with   supp $a\subset (x_0, x_1), (x_0>0)
 $. Let us   prove that $\frac{1}{q'}Ha$ is a $  (p,q,s)$-atom on ${\bf R}^+$. In fact,
 \par (i) we have already proved   supp $Ha\subset (x_0, x_1) $,
 \par (ii)  when $q=\infty $, noticing that $|Ha(x)|=\left|\frac {1}{x}\int^x_{x_0}a(t)dt\right|\leq  \|a\|_{L^\infty
 ({\bf R}^+)}\times \frac{x-x_0}{x}<\|a\|_{L^\infty  ({\bf R}^+)}$,
 we have
$$\|Ha\|_{L^\infty  ({\bf R}^+)}< \|a\|_{L^\infty  ({\bf R}^+)}\leq \left(  {x_1}-{x_0} \right)^{ -1/p};$$
 and when $1<q<\infty$, by the $L^q({\bf R^+})$ boundedness of $H$, we have
\begin{equation}
  \|\frac{1}{q'}Ha\|_{L^q({\bf R^+})}< \|a\|_{L^q({\bf R^+})}\leq \left(  {x_1}-{x_0} \right)^{1/q-1/p} ,
\end{equation}
\par (iii)~~ by the vanishing property of $a$, we have
\begin{eqnarray*}
\int^{+\infty}_{0}x^{\beta}Ha(x)dx
&=&\int^{x_1}_{x_0}x^{\beta}Ha(x)dx
\\
&=&\int^{x_1}_{x_0}x^{\beta -1}\int^{x}_{x_0}a(t)dtdx
\\
&=&\int^{x_1}_{x_0}a(t)
\int^{x_1}  _{t}x^{\beta -1}dxdt
\\
&=&
\left\{\begin{array}{cc}
\frac {1}{\beta}\int^{x_1}_{x_0}a(t)
\left(x_1^{\beta  }-t^{\beta}\right)dt, &  {\rm~~  if} ~   \beta =1, 2, ...,   s , \\
\int^{x_1}_{x_0}a(t)
\left(\ln x_1-\ln t \right)dt,& {\rm if}~~ \beta = 0
\end{array}
\right.
\\
&=& 0,{\rm~~  if} ~   \beta =0, 1, 2,  ...,  s  .
\end{eqnarray*}
Thus, combining (i), (ii) and (iii),  we see that $\frac{1}{q'}Ha$
is a $ (p,q,s)$-atom on ${\bf R}^+$. Let  $f \in LH ^ {p,q,s} ({\bf
R^+} )$, then
 $f=\sum _{j=1}^{\infty}\lambda _ja_j,$
where each $a_j$ is a  $ L-(p,q,s)$-atom on ${\bf R}^+$. So we have
$\frac{1}{q'}Hf=\sum _{j=1}^{\infty}\lambda _j\frac{1}{q'}Ha_j ,$
where  each $\frac{1}{q'}Ha_j$ is a  $  (p,q,s)$-atom on ${\bf
R}^+$. And by the definition: $\|\frac{1}{q'}Hf\|_{H ^{ p,q,s} ({\bf
R^+})} = \inf \left(\sum
 _{j=1}^{\infty}|\mu _j|^p\right)^{1/p}
,$ where the infimum is taken over all the decompositions
$\frac{1}{q'}Hf=\sum _{j=1}^{\infty}\mu  _jb_j,$ and  each $b_j$ is
a   $  (p,q,s)$-atom, we have
\begin{eqnarray*}
\|\frac{1}{q'}Hf\|_{H ^{p,q,s} ({\bf R^+})}
&=& \inf\limits
_{\stackrel{
\frac{1}{q'}Hf=\sum\limits_{j=1}^{\infty}\mu  _jb_j,}
{
{\rm each}
~ b_j
{\rm~ is~ a~  }
~ (p,q,s)-{\rm~ atom} }}
 \left(\sum
\limits_{j=1}^{\infty}|\mu _j|^p\right)^{1/p}
\\
&\leq & \inf
\limits
_{\stackrel{
\frac{1}{q'}Hf=\sum\limits_{j=1}^{\infty}\lambda _j\frac{1}{q'}Ha_j ,}
{{\rm each}~ a_j {\rm ~ is~ a~  }  L-(p,q,s)-{\rm ~ atom}}}
 \left(\sum
\limits_{j=1}^{\infty}|\lambda _j|^p\right)^{1/p}
\\
& &(
{\rm since~ each}~ \frac{1}{q'}Ha_j {\rm ~ is~ a~  }~
   (p,q,s)-{\rm ~ atom})
\\
&=& \inf \limits
_{\stackrel{
f=\sum\limits_{j=1}^{\infty}\lambda _ja_j,}
{{\rm each}~ a_j {\rm ~is~ a~  }~ L-  (p,q,s)-{\rm ~atom} }}
 \left(\sum
\limits_{j=1}^{\infty}|\lambda _j|^p\right)^{1/p}
\\
&=&\|f\|_{LH {p,q,s} ({\bf R^+})}.
\end{eqnarray*}
 Then
$$\|Hf\|_{H ^{p,q,s} ({\bf R^+})}\leq q'\|f\|^p_{LH ^{p,q,s} ({\bf R^+})}.$$
  \par  Thus, we finish the proof of Theorem 3.
\\
\par {\bf Proof of Theorem 4}~~
Let $a$ be a   $ (p,q,s)_{x^p}$-atom on ${\bf R}^+$ with supp
$a\subset (x_0, x_1)\subset  {\bf R}^+, (x_0>0) $. Let us prove that
\begin{eqnarray}
  \begin{array}{cc}
(1+p)^{-1/p}H^*a(x), &{\rm if }~~ q=\infty,\\
 (1-p)(1+p)^{1-1/p} H^*a(x), &{\rm if }~~ q=1, p\not= 1,\\
{\rm and }~\frac{1}{q}H^*a(x), &{\rm if }~~ 1<q<\infty, 0<p\leq 1,
\end{array}
\end{eqnarray}
are  the   $(p ,q,s-1)$-atoms
.
\par (i) Clearly,
  ${\rm supp} H^*a\subseteq (x_0, x_1)$.
  \par (ii)  When $q=\infty,$ and
  $x\in (x_0, x_1)$, noticing that
 $\|a\| _{L^\infty_{t^p}({\bf R}^+)}=\|a\| _{L^\infty ({\bf R}^+)},
$ and by (10), we have
\begin{eqnarray*}
 |H^*a(x)|
&\leq &  (x_1-x)\|a\| _{L^\infty  ({\bf R}^+)}
\\
&= &  (x_1-x )\|a\| _{L^\infty_{t^p} ({\bf R}^+)}
\\
&< & (x_1-x_0)
\left(\int^{x_1}_{x_0}  x^pdx \right)^{-1/p}
  \\
&= &
   (1+p)^{1/p}\frac{(x_1-x_0)
}{\left(x_1^{p+1}-x_0^{p+1}\right)^{1/p}}
\\
&< &
   (1+p)^{1/p}\frac{(x_1-x_0)
}{\left((x_1 -x_0)^{p+1}\right)^{1/p}}
 \end{eqnarray*}
\begin{equation}
~~~~~~~~~=~~
 (1+p)^{1/p} (x_1 -x_0) ^{-1/p} .
 \end{equation}
 When $  q =1, 0<p<1$,  by (10) and (11), we have,
\begin{eqnarray*}
\int^{\infty}_0|H^*a(x)|dx
   &\leq &
\int^{x_1}_{x_0}\left| \int^{x_1 }_{x} |a(t)|t^p \frac{1}{t^p}dt\right|dx
 \\
 &\leq &
\|a\| _{L^1_{t^p}({\bf R}^+)}
  \int^{x_1}_{x_0}    \frac{1}{x^{p }}
  dx
\\
 &\leq &
 \left(\int^{x_1 }_{x_0}  t^{p }dt
\right)^{1-1/p}
 \int^{x_1}_{x_0}    \frac{1}{x^{p }}
  dx
\\
&= &
\frac{1}{(1-p)(1+p)^{1-1/p}} \frac{  x_1 ^{1-p }-x_0^{1-p }  }
{\left(x_1^{1+p}-x_0^{1+p}\right)^{1/p-1}}
\\
&< &
\frac{1}{(1-p)(1+p)^{1-1/p}} \frac{ (x_1  -x_0)^{1-p} }{ (x_1 -x_0)^{(1+p)(1-p)/p} }
\\
&= &
\frac{1}{(1-p)(1+p)^{1-1/p}} (x_1  -x_0)^{1-1/p}  .
\end{eqnarray*}
  When $1< q <\infty,$ by (2), we have
\begin{equation}
\|\frac{1}{q}H^*a\|_{L^q({\bf R^+})}< \|a\|_{L^q({\bf R^+},t^q)}\leq \left(  {x_1}-{x_0} \right)^{1/q-1/p} .
\end{equation}
 ( iii)~~ By the vanishing property of $a$,  we have
\begin{eqnarray*}
\int^{+\infty}_{0}x^{\beta}H^*a(x)dx
&=&\int^{x_1}_{x_0}x^{\beta}H^*a(x)dx
 \\
&=&\int^{x_1}_{x_0}x^{\beta }\int^{x_1} _xa(t)dtdx
\\
&=&\int^{x_1}_{x_0}a(t)
\int^{t} _{x_0}x^{\beta }dxdt
\\
&=&
\frac {1}{\beta+1}\int^{x_1}_{x_0}a(t)
t^{\beta+1}dt
\\
&=& 0,{\rm ~~~~if}~~  \beta =0,1,2, ...,  s-1,
\end{eqnarray*}
when $0\leq \beta \leq s-1$. Thus, (21) has been proved. Therefore,
Theorem 4 follows from this   by analogous arguments to those in the
proof of Theorem 3.
 \par \section*
 {\bf  3.  Remarks}
  \par~~~~ {\bf Remark 1}~~ If the functions $f$ in Theorems 3 and   4 are finite linear combinations of
  corresponding atoms, i.e. $f=\sum _{j=1}^{k_0}\lambda _ja_j,$ then
     $${\rm the}~ "\leq "~{\rm    in}~ (8),(9) {\rm~will~ be~ changed~ to}~"< " .$$
  In fact, for Theorem 3,   suppose that there is $k_0<\infty$ such that
$$f=\sum\limits _{j=1}^{k_0}\lambda _ja_j,~
{\rm where~ each~} a_j ~{\rm is~ a}~   L-(p,q,s)-{\rm atom~ on}~{\bf R}^+,$$
 then, by the proof of theorem 3, for each $a_j, j=1,2,...,k_0$, we have
 $$\| Ha_j\|_{L^q({\bf R^+})}< q' \left(  {x_1}-{x_0} \right)^{1/q-1/p} ,$$
and it follows that there is $0<\epsilon _j<1 $ such that
 $$\| Ha_j\|_{L^q({\bf R^+})}< (q'-\epsilon _j) \left(  {x_1}-{x_0} \right)^{1/q-1/p} .$$
Let $\epsilon ={\rm min} \{\epsilon _1,\epsilon _2,...,\epsilon
_{k_0}\} $, then
$$\| Ha_j\|_{L^q({\bf R^+})}< (q'-\epsilon  ) \left(  {x_1}-{x_0} \right)^{1/q-1/p} .$$
As the arguments  in the proof of Theorem 3, we have that each
$\frac{1}{q'-\epsilon}Ha_j$ is a $  (p,q,s)$-atom on ${\bf R}^+$,
and
$$\|Hf\|_{H ^{p,q,s} ({\bf R^+})}\leq (q'-\epsilon)\|f\|^p_{LH ^{p,q,s} ({\bf R^+})}<q' \|f\|^p_{LH ^{p,q,s} ({\bf R^+})}.$$
\par Similar  arguments   above are suitable for the case of  Theorem 4.
\par {\bf Remark 2}~~ If we drop  the restriction $x_0>0$ in the definitions of Hardy spaces (in Definition 1), then
$${\rm the}~ "<"~{\rm    in}~ (6),(7) {\rm~will~ be~ changed~ to}~"\leq " .$$
This is because of that   the "$<$ " in (12) and (19) will be
changed to "$ \leq $" if $x_0=0 $ in the proofs of Theorems 1  and
  2.
\par {\bf Remark 3}~~ If the functions $f$ in Theorems 3 and   4 are finite linear combinations of
corresponding atoms, i.e. $f=\sum _{j=1}^{k_0}\lambda _ja_j,$ then,
even if   dropping the restriction $x_0>0$ in the definitions of
Hardy spaces (in Definition 1), we have that:
     $${\rm the}~ "\leq "~{\rm    in}~ (8)  {\rm~will~ be~ changed~ to}~"< " $$
 when $1< q< \infty,   $ and
  $${\rm the}~ "\leq "~{\rm    in}~  (9) {\rm~will~ be~ changed~ to}~"< " $$
when $ 1< q\leq \infty  $. These follow from   the analogous
arguments to those of Remark 1, since  the "$<$" in (20),(23),
 and the first "$<$" in (22) hold still when $x_0=0$.

 \par Shunchao  Long
\par Mathematics Department,
\par Xiangtan University,
\par Xiangtan , 411105,  China
\par E-mail address: sclong@xtu.edu.cn
\end{document}